\documentclass[12pt]{amsart}
\usepackage{amsfonts}
\usepackage{amssymb}
\numberwithin{equation}{section}
\theoremstyle{plain}
 \newtheorem{thm}{Theorem}[section]
 \newtheorem{lem}{Lemma}[section]
 \newtheorem{cor}{Corollary}[section] 
 \newtheorem{prop}{Proposition}[section]
\theoremstyle{definition}
 \newtheorem{defn}{Definition}[section]
  \newtheorem{rem}{Remark}[section] 

\newcommand{\al}{\alpha}

\newcommand{\ep}{\varepsilon}

\newcommand{\ch}{\widehat \mu}
\newcommand{\q}{\quad}

\newcommand{\eqd}{\overset{\mathrm d}{=}}

\newcommand{\G}{\mathbb{G}}
\newcommand{\R}{\mathbb{R}}
\newcommand{\N}{\mathbb{N}}
\newcommand{\Z}{\mathbb{Z}}

\newcommand{\D}{\mathbb{D}}

\newcommand{\les}{\leqslant}
\newcommand{\ges}{\geqslant}

\newcommand{\cl}{\colon}

\begin{document}

\footnote[0]{T. Watanabe: Center for Mathematical Sciences The Univ. of Aizu. Aizu-Wakamatsu 965-8580, Japan. \\
e-mail: t-watanb@u-aizu.ac.jp}

\title{ Escape rates for multi-dimensional  shift selfsimilar additive sequences  }
\author{Toshiro Watanabe}
\maketitle

{\small {\bf Abstract }
First the relation between shift selfsimilar additive  sequences and  stationary  sequences of Ornstein-Uhlenbeck type (OU type) on $\R^d$ is shown and then the rates of escape for shift selfsimilar additive  sequences are discussed. As a corollary, fundamental problems on recurrence of stationary  sequences of OU type are solved. Some applications to laws of the iterated logarithm for strictly stable L\'evy 
processes on $\R^d$ and independent Brownian motions are given. }\\

\medskip

\bigskip


\section{Introduction and main results}
In what follows,  denote the Euclidean inner product of $z$ and $x$ in $\R^d$ and the Euclidean norm of $x$ by $\langle z,x \rangle$ and $|x|$, respectively. Let $\Z:=\{0,\pm 1, \pm 2, \dots \} $, $\Z_+:=\{0, 1,  2, \dots \} $,  $\N:= \{ 1,  2, 3, \dots \} $, and $\R_+:=[0,\infty)$.  The symbol $\delta_a(dx)$ stands for the delta measure at $a$ in $\R^d$.  For an $\R^d$-valued random variable $X$, let  $S_{X}$ be the support of the distribution of $X$. Denote 
the characteristic function of a probability distribution $\mu$ on $\R^d$  by $\ch(z)$, namely,
\begin{equation}
 \ch(z): = \int_{\R^d}e^{i\langle z,x \rangle}\mu(dx).\nonumber
\end{equation} 
 We use the words ``increase" and ``decrease" in the wide sense allowing flatness. For positive  functions $f(x)$ and $g(x)$ on $(0,\ep)$ for some $\ep \in (0,1)$, we define the relation $f(x) \sim g(x)$ as $x \to 0+$ by $\lim_{x \to 0+}f(x)/g(x) =1$ and the relation $f(x) \asymp g(x)$ by $0 <\liminf_{x \to 0+}f(x)/g(x) \les \limsup_{x \to 0+}f(x)/g(x) < \infty$.

\begin{defn} Let $a >1$.  A random sequence $\{W(n),n \in \Z \}$ on $\R^d$ is called a {\it shift $a$-selfsimilar additive  sequence } if the following  hold.

(1) The sequence  $\{W(n)\}$ is shift $a$-selfsimilar, that is, 
$$ \{W(n+1), n \in \Z
 \}\stackrel{\rm d}{ =} \{aW(n), n \in \Z \}, $$
where $\stackrel{\rm d}{ =}$ denotes the equality in  finite-dimensional distributions.

(2) The sequence  $\{W(n)\}$ has independent increments, that is, for every $n \in \Z $, $\{W(k), k \les n\}$ and  $W(n+1)-W(n)$ are independent. 
\end{defn}

\begin{defn} Let $0 < b < 1$. 
A probability distribution $\mu$ on $\R^d$ is said to be {\em $b$-decomposable } if there exists a   probability  distribution $\rho$ on $\R^d$ such that
\begin{equation}
 \ch(z) = \ch(bz)\widehat \rho(z).
\end{equation}  
\end{defn}
\begin{defn}
A probability distribution $\zeta$ on $\R^d$ is said to be {\it full} if its support is not contained in any  hyperplane in $\R^d$.
\end{defn}

 The decomposition (1.1) is equivalent to 
$$ \widehat {\mu}(z) = \prod_{n=0}^{\infty}{\widehat \rho}(b^nz). $$
The infinite product above converges if and only if $\rho$ has  finite log-moment, namely,
\begin{equation}
 \int_{\R^d} \log(2+|x|)\rho(dx) < \infty. 
\end{equation} 
See  Lemma 1 of Bunge \cite{b97}. The distribution $\mu$ is uniquely determined by $\rho$, but $\rho$ is not always uniquely determined by $\mu$.  If  the support of $\rho$ is a compact set, then $\rho$ is  uniquely determined by $\mu$.

\begin{rem} The author characterized the shift $a$-selfsimilar additive  sequence in Theorem 2.1 of \cite{w00a} as follows.

 (i)  Let $\{W(n), n \in \Z \}$ be an $\R^d$-valued shift $a$-selfsimilar additive sequence. Let $\mu$ and $\rho$ be the distributions of $W(0)$ and $W(0)-W(-1)$, respectively.
Then $\mu$ is a $b$-decomposable distribution satisfying (1.1) with $b=a^{-1}$.

(ii)  Conversely, if a probability distribution $\rho$ with finite log-moment  is given, then there is a 
unique $($in law$)$ shift $a$-selfsimilar  additive  sequence $\{W(n), n \in \Z \}$ such that  $\rho$ is the distribution of  $W(0)-W(-1)$.  
That is, for every $a^{-1}$-decomposable distribution $\mu$ on $\R^d$, there is a 
 $($not necessarily unique in law$)$ shift $a$-selfsimilar  additive  sequence $\{W(n), n \in \Z \}$ such that $\mu$ is the distribution of  $W(0)$.  
\end{rem}
 
\begin{defn} 
 Let $a >1$.  A random  sequence $\{Y(n),n \in \Z \}$ on $\R^d$ is called a {\it stationary  sequence of Ornstein-Uhlenbeck type} (OU type in short) with parameter $a$ if  the following hold.

(1) The sequence $\{Y(n)\}$ is  stationary, that is, 
$$ \{Y(n+1), n \in \Z
 \}\stackrel{\rm d}{ =} \{Y(n), n \in \Z \}. $$

(2) For every $n \in \Z $, $\{Y(k), k \les n\}$ and  $Y(n+1)-a^{-1}Y(n)$ are independent. 
\end{defn}

\begin{defn}  Let $\{Y(n), n \in \Z \}$ be a stationary  sequence of OU type  on $\R^d$ with parameter $a$. A point $x \in \R^d$ is a {\it recurrent point} for $\{Y(n)\}$ if 
\begin{equation}
\liminf_{n \to   \infty}|Y(n) -x| = 0 \mbox{ a.s..} \nonumber
\end{equation}
We define  a class $\G$ as the totality of positive decreasing sequences $g(n)$ on  $\Z_+ $. A recurrent  point $x$ is of {\it type A} for $\{Y(n)\}$ if there exists $g \in \G$ such that
\begin{equation}
\liminf_{n \to   \infty}\frac{|Y(n)-x|}{g(n)} = 1 \mbox{ a.s.}.
\end{equation}
A recurrent  point $x$ is of {\it type B} for $\{Y(n)\}$ if there  exists no $g \in \G$ satisfying (1.3).

\end{defn}
 
\begin{defn} 
Let $f(r)$ be a nonnegative
increasing function on $(0,1)$.  The function $f(r)$ is called of {\it dominated variation } as $r \to  0+$ if it is positive on $(0,\ep)$ for some $\ep \in (0,1)$ and if  $\limsup_{r \to 0+} f(2r)/f(r) < \infty,$ that is, $f(2r) \asymp f(r)$ as $ r \to 0+$.
\end{defn}

Throughout the paper, let $a >1$ and  we assume that the  distributions of  $W(0)$ and $Y(0)$ are full on $\R^d$. Motivated by the results of  Maejima and Sato \cite{ms} and Pruitt \cite{p90}, the author \cite{w02a} introduced the shift selfsimilar additive  sequences and started the research on the asymptotic behaviors of those sequences. Main problems on limsup behaviors for those sequences were solved by \cite{w02a}.  Following this, he gave applications to laws of the iterated logarithm for Brownian motions on nested fractals  and to  determining the  exact Hausdorff and packing measures of random fractals on Galton-Watson trees in a series of papers \cite{w02b,w04,w07}. In this paper, we achieve the study on the rates of escape (liminf behaviors) for multi-dimensional shift selfsimilar additive  sequences. At the same time, basic problems on recurrence of  stationary  sequences of OU type are discussed.  Our main results are as follows. 

\begin{thm}  Let  $\{W(n)\}$ be a  shift $a$-selfsimilar additive  sequence on $\R^d$. 

{\em (i)} Let $D:= \inf_{y \in S_{W(0)}}|y|$.   Then we have
\begin{equation}
\liminf_{n \to   \infty}\frac{|W(n)|}{a^n} = D \mbox{ {\em a.s..}}\nonumber
\end{equation}

{\em (ii)} There exists $g \in \G$ such that
\begin{equation}
\liminf_{n \to   \infty}\frac{|W(n)|}{a^ng(n)} = 1 \mbox{ {\em a.s.}}
\end{equation}
if and only if  $P(|W(0)| \les r)$ is not of dominated variation as $r \to 0+$.
\end{thm}

\begin{cor} Let $\{Y(n)\}$ be a stationary  sequence of OU type on $\R^d$ with parameter $a$. 

{\em (i)}  Let $D(x,S_{Y(0)}) := \inf_{y \in S_{Y(0)}}|x-y|$. Then we have
\begin{equation}
\liminf_{n \to   \infty}|Y(n) -x| = D(x,S_{Y(0)}) \mbox{ {\em a.s..}}\nonumber
\end{equation}
 Thus a point $x \in \R^d$ is a recurrent point for $\{Y(n)\}$ if and only if $x \in S_{Y(0)}.$

{\em (ii)}   A recurrent point $x$ is of type A for $\{Y(n)\}$ if and only if  $P(|Y(0)-x| \les r)$ is not of dominated variation as $r \to 0+$.

\end{cor}

\begin{prop}  Let  $\{W(n)\}$ be a  shift $a$-selfsimilar additive  sequence on $\R_+^d$. Then there exists $g \in \G$ satisfying $(1.4)$ if and only if  $P(|W(0)-W(-1)|=0)=0$.

\end{prop}

\begin{thm}  Let $\{W(n)\}$ be a  shift $a$-selfsimilar additive  sequence on $\R^d$. Suppose that  $F(r) :=P(|W(0)| \les r)$ is of dominated variation as $r \to 0+$. Let $g \in  \G $.  If
\begin{equation}
 \sum_{n=0}^{\infty}F(g(n)) = \infty \q ( \mbox{{\em resp.}} \q < \infty ),\nonumber
\end{equation}
then
\begin{equation} 
\liminf_{n \to  \infty}\frac{|W(n)|}{a^ng(n)} = 0 \q ( \mbox{{\em resp.}}  \q = \infty )\q \mbox{ {\em a.s..}}
\end{equation}
\end{thm}

\begin{cor}  Let  $\{Y(n)\}$ be a stationary  sequence of OU type on $\R^d$ with parameter $a$. Suppose that $x \in S_{Y(0)}$ is of type B for $\{Y(n)\}$ and hence $F_x(r) :=P(|Y(0)-x| \les r)$ is of dominated variation as $r \to 0+$. Let $g \in  \G $.  If
\begin{equation}
 \sum_{n=0}^{\infty}F_x(g(n)) = \infty \q ( \mbox{{\em  resp.}} \q < \infty ),\nonumber
\end{equation}
then
\begin{equation}
 \liminf_{n \to  \infty}\frac{|Y(n)-x|}{g(n)} = 0 \q ( \mbox{{\em resp.}}  \q = \infty )\q \mbox{ {\em a.s..}}
\end{equation}
\end{cor}
\medskip
Under the assumption that $\lambda:= P(|W(0)-W(-1)|=0)>0$,
define a Laplace transform $ L_W(u)$ for $ u \ges 0$ as
$$L_W(u):= E(\exp(-u|W(0)-W(-1)|))$$ 
and a regularly varying function $K_{W}(r)$ on $(0,1]$ with index $-\log \lambda /\log a$  as $r \to 0+$  as
$$ K_{W}(r):=  r^{-\log \lambda/ \log a} \exp( \int_1^{r^{-1}} \frac{\log L_W(u)-\log \lambda }{u \log a} du ). $$

\begin{prop}  Let $\{W(n)\}$ be a  shift $a$-selfsimilar additive  sequence on $\R_+^d$. Assume that $P(|W(0)-W(-1)|=0)> 0$. 
Let $g \in \G $.  If
$$ \sum_{n=0}^{\infty}K_{W}(g(n)\wedge 1) = \infty \quad ( \mbox{{\em resp.}} \quad < \infty ),$$
then $(1.5)$ holds.

\end{prop}

In Sect.\ 2, we discuss the series representation of $\{W(n)\}$ and $\{Y(n)\}$. In Sect.\ 3, we prove the main results stated in Sect.\ 1. In Sect.\ 4 and Sect.\ 5,  we give some applications of main results to laws of the iterated logarithm for strictly stable L\'evy processes and  independent Brownian motions, respectively. We shall discuss the rate of access of $\{Y(n)\}$ for recurrent points in a forthcoming paper.

\section{Series representation of $\{W(n)\}$ and $\{Y(n)\}$ }
Following Sato \cite{s91,s99}, we define a  Sato process. Moreover, we add the definition of a  process of OU type. 
\medskip
\begin{defn} Let $H >0$. A stochastic process $\{X(t), t \in \R_+\}$ is called a {\it Sato process} with exponent $H$ if the following hold.

(1) $\{X(t), t \in \R_+\}$ is {\it $H$-selfsimilar} with exponent $H$, that is, for every $c >0$
 $$\{X(ct), t \in \R_+\} \stackrel{\rm d}{ =}  \{c^HX(t), t \in \R_+\}.$$ 

(2) $\{X(t), t \in \R_+\}$ has independent increments.

(3) $\{X(t), t \in \R_+\}$ is stochastically continuous with c\`adl\`ag paths. 
\end{defn}
\begin{defn} Let $\{Z(t), t \in \R_+\}$ be a L\'evy process on $\R^d$. 
A stochastic process $\{U(t), t \in \R_+\}$ is called a  {\it  process of OU type} with initial state $U(0)$ with  parameter $c > 0$ if it is the solution of 
$$U(t) = U(0) +Z(t) -c\int_0^t U(s)ds,$$
which is 
$$U(t) = e^{-ct}(U(0)+\int_0^t e^{cs}dZ(s)).$$

\end{defn}

Let $b > 1$. Then the sequence $\{X(b^n), n \in \Z\}$ for a Sato process $\{X(t), t \in \R_+\}$  with exponent $H$ is a shift $b^H$-selfsimilar additive  sequence. Every stationary process $\{U(t)\}$ of OU type for  $t \in \R_+$  with parameter $c$ can be extended to the stationary process $\{U(t)\}$ for $t \in \R^1.$ Its periodic observation  $\{U(n), n \in \Z \}$  is a stationary  sequence of OU type  with parameter $e^c$.  There is a one to one correspondence between  stationary  processes of OU type  and Sato processes through the Lamperti transformation.  Moreover, those two processes are represented by integrals of L\'evy processes. Refer to Jeanblanc et  al. \cite{jpy} and also to Maejima and Sato \cite{ms03}. In this section, we give discrete time analogue of those facts. 

\begin{prop}   There is a one to one correspondence between a stationary  sequence $\{Y(n)\}$ of OU type  with parameter $a$ and a shift $a$-selfsimilar additive  sequence $\{W(n)\}$ through the Lamperti transformation :
\begin{equation}
W(n)=a^n Y(n).
\end{equation}
 \end{prop}

{\it Proof } Let $\{W(n)\}$ and $\{Y(n)\}$ be two random sequences satisfying the Lamperti transformation (2.1). Then  $\{W(n)\}$ has shift $a$-selfsimilarity if and only if $\{Y(n)\}$ has stationarity. Moreover, for every $n \in \Z $, $\{W(k), k \les n\}$ and  $W(n+1)-W(n)$ are independent if and only if, for every $n \in \Z $,  $\{Y(k), k \les n\}$ and  $Y(n+1)-a^{-1}Y(n)$ are independent. Thus the theorem has been proved. \qed

\begin{rem} By virtue of the above proposition and Remark 1.1,  the stationary  sequence of OU type with parameter $a$ is characterized as follows.

 (i) Let $\{Y(n), n \in \Z \}$ be an $\R^d$-valued stationary  sequence of OU type with parameter $a$. Let $\mu$ and $\rho$ be the distributions of $Y(0)$ and $Y(0)-a^{-1}Y(-1)$, respectively.
Then $\mu$ is a $b$-decomposable distribution satisfying (1.1)  with $b=a^{-1}$.

 (ii) Conversely, if a probability distribution $\rho$ with finite log-moment  is given, then there is a 
unique $($in law$)$ stationary  sequence $\{Y(n), n \in \Z \}$ of OU type with parameter $a$ such that  $\rho$ is the distribution of  $Y(0)-a^{-1}Y(-1)$.  
That is, for every $a^{-1}$-decomposable distribution $\mu$ on $\R^d$, there is a 
 $($not necessarily unique in law$)$ stationary  sequence $\{Y(n), n \in \Z \}$ of OU type with parameter $a$ such that $\mu$ is the distribution of  $Y(0)$.  
\end{rem}

The following lemma is due to (i) of Theorem 2.1 of Watanabe \cite{w02a}.
\begin{lem} Let $\{W(n)\}$ be a shift $a$-selfsimilar additive  sequence on $\R^d$. Then we have
\begin{equation}
\lim_{n \to   -\infty}|W(n)| = 0 \mbox{ {\em a.s.}}.\nonumber
\end{equation}
\end{lem}

The following lemma is due to Lemma 1 of Bunge \cite{b97}.
\begin{lem}  Let $\{X_n\}_{n=-\infty}^{\infty}$ be $\R^d$-valued i.i.d. ramdom variables. Then the sum $\sum_{j=-\infty}^0 a^j X_j $ is convergent almost surely if and only if $X_0$ has  finite log-moment, that is,  
\begin{equation}
 E(\log(|X_0| +2)) < \infty. \nonumber
\end{equation} 
\end{lem}

\begin{thm}   Let $\{W(n)\}$ be a shift $a$-selfsimilar additive  sequence on $\R^d$. Define random variables $\{X_n\}_{n=-\infty}^{\infty}$  as 
$$X_n := a^{-n}(W(n)-W(n-1)).$$ 
Then  $\{X_n\}$ are $\R^d$-valued i.i.d. ramdom variables with $X_0$ having  finite log-moment.
The sequence $\{W(n)\}$ is expressed almost surely as
\begin{equation}
W(n) = \sum_{j=-\infty}^n a^j X_j.  
\end{equation}
Conversely, let $\{X_n\}_{n=-\infty}^{\infty}$ be an $\R^d$-valued i.i.d. ramdom variables with $X_0$ having  finite log-moment. Then the sequence $\{W(n)\}$ defined by $(2.2)$  is  a shift $a$-selfsimilar additive  sequence.
 \end{thm}

\medskip

{\it Proof } Let $X_n := a^{-n}(W(n)-W(n-1))$ for $n \in \Z$. Then, by the shift $a$-selfsimilarity of $\{W(n)\}$, we see that $\{X_n\}$ are identically distributed. Moreover, by the
 independence of increments, $\{X_n\}$ are independent. Thus $\{X_n\}$ are i.i.d.. We have, for $k \ges 1$,
\begin{equation}
W(n)-W(n-k) = \sum_{j=n-k+1}^{n} a^jX_j.\nonumber
\end{equation}
Thus we obtain (2.2) from Lemmas 2.1 and 2.2 and see that
 $X_0$ has finite log-moment. The converse statement is obviously true.   \qed

\medskip
 
We obtain the following corollary directly from Theorem 2.1 and Proposition 2.1, and thereby  its proof is omitted.
 
\begin{cor}  Let  $\{Y(n)\}$ be a stationary  sequence of OU type with parameter $a$ on $\R^d$. Define random variables $\{X_n\}_{n=-\infty}^{\infty}$  as 
$$X_n : = Y(n)-a^{-1}Y(n-1).$$ 
Then  $\{X_n\}$ are $\R^d$-valued i.i.d. ramdom variables with $X_0$ having  finite log-moment and the sequence $\{Y(n)\}$  is expressed almost surely  as
\begin{equation}
Y(n) = \sum_{j=-\infty}^n a^{j-n} X_{j}. 
\end{equation}
Conversely, let $\{X_n\}_{n=-\infty}^{\infty}$ be an $\R^d$-valued i.i.d. ramdom variables with $X_0$ having  finite log-moment. Then  the sequence $\{Y(n)\}$ defined by $(2.3)$ is a stationary  sequence of OU type with parameter $a$. 
\end{cor}

\begin{prop} Every stationary sequence $\{Y(n)\}$ of OU type is ergodic.  It holds that, for every  $\delta >0$,
\begin{equation}
\lim_{n \to  \infty}\frac{1}{n}\sum_{k=1}^{n}1_{\{|Y(k)-x| \les \delta\}} = P(|Y(0)-x| \les \delta)    \mbox{ {\em a.s.}}\nonumber
\end{equation}
and 
\begin{equation}
\lim_{n \to  -\infty}\frac{1}{|n|}\sum_{k=n}^{-1}1_{\{|Y(k)-x| \les \delta\}} = P(|Y(0)-x| \les \delta)\mbox{ {\em  a.s.}}.\nonumber
\end{equation}
\end{prop}
\medskip

{\it Proof } The first assertion follows from Propositions 3.1 and 3.3 of Breiman \cite{b68} and Corollary 2.1. The second assertion is clear by virtue of the ergodic theorem for stationary sequences. \qed

\section{Proofs of main results}
In this section, let $a >1$ and $\{W(n)\}$ be a shift $a$-selfsimilar additive  sequence on $\R^d$.  We establish the results on the rate of escape for $\{W(n)\}$.  Corolleries 1.1 and 1.2 for  stationary  sequences $\{Y(n)\}$ of OU type can be proved  by letting $W(n):= a^n(Y(n)-x)$ for $n \in \Z$. Thus their proofs are omitted.
We define  a class $\G_1$ of sequences $g(n)$ as 
 $$\G_1: = \{ g(n) : g(n) \in \G,   g(n+1)/g(n) \ges \gamma^{-1}  \mbox{ for some $\gamma \in (1,a)$  on } \Z_+ \}.$$ 
Let $b = a^{-1}$. Define  functions $ H_b(u)$ for $u \ges 0$ and $I(c)$ for $c > 0$ as
\begin{equation} 
 H_b(u) : =  \sum_{k=0}^{\infty}(1-  \cos(b^k u) )\nonumber
\end{equation} 
and
$$ I(c): = \sum_{n=0}^{\infty}\int_0^{2\pi} \exp(c\sum_{k=0}^n( \cos(a^{k}u) - 1)) du.$$
\begin{lem}  Let $g \in \G$. If $\sum_{n=0}^{\infty}P(|W(0)| \les g(n)) < \infty$, then 
\begin{equation}
\liminf_{n \to   \infty}\frac{|W(n)|}{a^ng(n)} \ges 1 \mbox{ {\em  a.s.}}.
\end{equation}
\end{lem}
{\it Proof } Suppose that $\sum_{n=0}^{\infty}P(|W(0)| \les g(n)) < \infty$.  Then we have by the shift selfsimilarity
\begin{equation}
\sum_{n=0}^{\infty}P(|W(n)| \les a^ng(n)) < \infty.\nonumber
\end{equation}
By virtue of the Borel-Cantelli lemma, we have (3.1). \qed

\medskip
The following lemma is due to (i) of Lemma 3.3 of Watanabe \cite{w02a}.
\begin{lem}  We have  $I(c) < \infty$ for every $c > 0$. 
\end{lem}
\begin{lem}   We have for every  $\delta >0$ 
\begin{equation}
\sum_{n=0}^{\infty}P(|W(n)| \les \delta) < \infty.
\end{equation}
\end{lem}
{\it Proof } Let $\mu$ and $\rho$ be the distributions of $W(0)$ and $W(0)-W(-1)$, respectively. Let $\bar{\rho}$ be the reflection of $\rho$, that is, $\bar{\rho}(dx):=\rho(-dx)$ and $\rho*\bar{\rho}$ be the convolution of $\rho$ and $\bar{\rho}$.   We find from  Lemma 5 of \cite{w00a} that $\rho$ and also $\rho*\bar{\rho}$ are full on $\R^d$. Thus there exist disjoint  closed balls $D_j$ for $1 \les j \les d$ in $\R^d$ such that $C_j :=\rho*\bar{\rho}(D_j) > 0$ and, for any choice of $x_j \in D_j,$  $\{x_j\}_{j=1}^d$ is a basis of $\R^d$.  Let $x_j ={}^t(x_{1j},x_{2j},\ldots,x_{dj})$ and define a real $d\times d$ matrix $X$ as $X = ( x_{ij} ) $. Clearly, if $x_j \in D_j$ for $1 \les j \les d$, then
 $ \det X \neq 0 $. Choose $a_0 >0$ such that $\sup \{ a_0 |x_j| : x_j \in D_j \mbox{ for }1 \les j \les d \} < 2\pi$.  Let $b=a^{-1}$. 
 Note that $\int_{\R^d} \cos \langle z,x\rangle \rho*\bar{\rho}(dx)= |\widehat {\rho}(z)|^2 \ges 0$.
Using the inequality 
$x \les e^{x-1}$ for $x\ges 0$, we have
\begin{equation} 
\begin{split}
 &  |\widehat {\mu}(z)|^2 = \prod_{k=0}^{\infty} \int_{\R^d} \cos \langle b^kz,x\rangle \rho*\bar{\rho}(dx) \nonumber \\ 
 \les &\exp ( \sum_{k=0}^{\infty} \int_{\R^d} (\cos \langle b^kz,x\rangle - 1 ) \rho*\bar{\rho}(dx))\nonumber \\
  = &\exp (- \int_{\R^d} H_b(|\langle z,x \rangle|) \rho*\bar{\rho}(dx) ).\nonumber
\end{split} 
\end{equation} 
Thus, applying Jensen's inequality and changing variables as $\langle z,x_j \rangle = u_j$ for $1 \leq j \leq d$, we obtain    that
\begin{equation} 
\begin{split}
&\sum_{n=0}^{\infty} \int_{|z| \les a_0}|\widehat {\mu}(a^{n}z)| dz \\ 
& \les \sum_{n=0}^{\infty} \int_{|z| \les a_0}\exp\{- 2^{-1} \int_{\R^d}H_b(|\langle a^{n} z,x \rangle|)\rho*\bar{\rho}(dx)\}dz  \\
& \les \sum_{n=0}^{\infty} \int_{|z| \les a_0}\exp\{- 2^{-1} \sum_{j=1}^d\int_{D_j}H_b(|\langle a^{n} z,x_j \rangle|)\rho*\bar{\rho}(dx_j)\}dz  \\
& \les \sum_{n=0}^{\infty} \int_{|z| \les a_0}dz \prod_{j=1}^d C_j^{-1}\int_{D_j}\exp\{ -2^{-1} C_j H_b(|\langle  a^{n} z,x_j \rangle|) \}\rho*\bar{\rho}(dx_j) \\ 
& \les (\prod_{j=1}^d C_j^{-1})(4\pi)^{d-1}2I(2^{-1}C_1)\int_{D_1\times\cdots\times D_d}|\det X |^{-1}\prod_{j=1}^d  \rho*\bar{\rho}(dx_j) < \infty. 
\end{split}
\end{equation}
In the last inequality  above,  we used Lemma 3.2.  Let $y =(y_j)_{j=1}^d \in \R^d.$ Define a function $f_c(y)$ on $\R^d$ for $c > 0$ as
$$ f_c(y) = \prod_{j=1}^d \left(\frac{2 \sin(cy_j/2)}{cy_j}\right)^2$$
with the understanding that $(\sin 0)/0 = 1$. Then the Fourier transform $\widehat f_c(z)$ of 
$f_c(y)$ is given by
$$\widehat f_c(z) = \int_{\R^d} \exp(i\langle z, y \rangle) f_c(y) dy$$
$$      = (2\pi c^{-1})^{d} \prod_{j=1}^d (1 - c^{-1}|z_j|) 1_{[-c,c]}(z_j),$$
where $1_{[-c,c]}(x)$ is the indicator function of the interval $[-c,c]$.
We see from the Perseval's equality that
\begin{equation} 
\begin{split}
 \sum_{n=0}^{\infty} E( f_c(W(n))) &= (2 \pi)^{-d} \sum_{n=0}^{\infty}\int_{\R^d}\widehat \mu(a^nz) \widehat f_c(-z)dz \\ \nonumber
 &\les  c^{-d}\sum_{n=0}^{\infty}\int_{|z| \leq \sqrt d c }|\widehat \mu(a^nz)| dz.\nonumber
\end{split}
\end{equation}
Hence we find that (3.3) implies that $ \sum_{n=0}^{\infty} E( f_c(W(n))) < \infty $ for any $c \in (0,a_0/\sqrt d).$ 
Thus the lemma is true.    \qed

\begin{lem}   We have for every $\delta >0$ and $c > 1$
\begin{equation}
\sum_{n=0}^{\infty}P(|W(0)| \les \delta c^{-n}) < \infty.
\end{equation}
\end{lem}
{\it Proof } 
 We have by (3.2) of  Lemma 3.3 and the shift selfsimilarity
\begin{equation}
\sum_{n=0}^{\infty}P(|W(0)| \les \delta a^{-n}) < \infty.\nonumber
\end{equation}
Thus we see that for every  $\delta >0$
\begin{equation}
\int_0^1 P(|W(0)| \les \delta x) \frac{dx}{x} < \infty.\nonumber
\end{equation}
Hence we have (3.4) for every   $c > 1$ and $\delta >0$.  \qed

\begin{lem} Let $g \in \G$. There exits $g^* \in \G_1$ such that $g(n) \les g^*(n)$ on $\Z_+$ and, for every $\delta > 0$, 
\begin{equation}
\sum_{n \in J}P(|W(0)| \les \delta g^*(n)) < \infty,
\end{equation}
where $J:=\{n \in\Z_+ : g(n) < g^*(n)\}.$
\end{lem}
{\it Proof } Fix $\gamma \in (1,a)$.  We define $g^* \in \G_1$  inductively. First let $g^*(0)= g(0)$. Assume that $g^*(n)$ is defined for $0 \les n \les k$.
If $g(k+1) < \gamma^{-1} g^*(k)$, then define $g^*(k+1):= \gamma^{-1} g^*(k).$  If  $g(k+1) \ges \gamma^{-1} g^*(k)$, then define $g^*(k+1):= g(k+1)$.
Thus $J= \{n \in \N : g(n) < \gamma^{-1}g^*(n-1) \}$. Hence  we obtain from  Lemma 3.4 that 
\begin{equation}
\sum_{n \in J}P(|W(0)| \les g^*(n)) \les  \sum_{n =0}^{\infty}P(|W(0)| \les \gamma^{-n} g(0)) < \infty.\nonumber
\end{equation}
Thus (3.5) is true. \qed

\begin{lem}  Let  $g \in \G_1$. If $\sum_{n=0}^{\infty}P(|W(0)| \les g(n)) = \infty$, then 
\begin{equation}
\liminf_{n \to   \infty}\frac{|W(n)|}{a^ng(n)} \les 1 \mbox{ {\em  a.s.}}.
\end{equation}
\end{lem}
{\it Proof } 
Suppose that 
$$ \sum_{n=0}^{\infty} P( |W(0)| \les g(n) ) = \infty.$$
Then we get by the shift selfsimilarity
\begin{equation}
 \sum_{n=0}^{\infty} P( |W(n)| \les a^{n}g(n) ) = \infty. 
\end{equation}
Define  events $A_n$ and a sequence $p_n$ with $\ell \in  \N $ as
$$ A_n := \{ \omega : |W(n)| \les a^{n}g(n) \}$$
and
$$ p_n := P(|W(j\ell) -W(n\ell)| > a^{j\ell}g(j\ell) + a^{n\ell}g(n\ell)\quad \mbox{ for } \forall j \ges n+1 ).$$
We find from (3.7) that there is $j$ with $0 \les j \les \ell-1$ such that
\begin{equation}
\sum_{n=0}^{\infty} P(A_{n\ell+j}) = \infty.
\end{equation}
 Without loss of generality, we can assume that $j=0$. We have 
\begin{equation}
1 \ges P(\bigcup_{n=0}^{\infty}A_{n\ell}) \ges \sum_{n=0}^{\infty}P((\bigcup_{j=n+1}^{\infty}A_{j\ell})^c\cap A_{n\ell}) \ges \sum_{n=0}^{\infty}P(A_{n\ell})p_n.
\end{equation}
In the last inequality above, we used inequalities $|W(j\ell)| \ges |W(j\ell) -W(n\ell)| -|W(n\ell)|$ for $j \ges n+1 $ and then the independence of increments.
We see from the shift selfsimilarity that
$$p_n = P(|W(k\ell) -W(0)| > a^{k\ell}g((n+k)\ell) + g(n\ell)\quad \mbox{ for } \forall k \in \N ).$$
Since $g(n)$ is decreasing, we find that $p_n$ is increasing and hence see from (3.8)  and (3.9) that $p_n = 0$ for every $n \in \Z_+$.
 Choose $\ep >0 $ arbitrarily. Define events $B_n$ and $C_n$ and a sequence  $q_n$ for $n \in \Z_+$ as 
$$ B_n:= \{ \omega : |W(n\ell))-W(0)| \les  (1+\ep)a^{n\ell}g(n\ell) \},$$
$$ C_n:=\{ \omega : |W((n+k)\ell) -W(0) |> (1+\ep)a^{(n+k)\ell}g((n+k)\ell) \quad \mbox{ for } \forall k \in \N  \},$$
and
$$ q_n := P(B_n \cap C_n).$$
Let $a_{n,k} :=a^{k\ell}g((n+k)\ell)-g(n\ell)$ and $b_{n,k}:= a^{k\ell}g((n+k)\ell) + g(n\ell)$ for $k \in \N $ and $n \in \Z_+$. Note from the definition of the class $\G_1$ that $a^{\ell}g((n+1)\ell) /g(n\ell) \ges (a/\gamma)^{\ell}$ with some $\gamma \in (1,a)$ for 
$n \in \Z_+$.  For any $\ep >0$, we can choose $\ell \in \N$ sufficiently large such that 
$(1+\ep)a_{n,k}   \ges b_{n,k}$ for $k \in \N $ and $n \in \Z_+$.
By using inequalities
$$ |W((n+k)\ell)-W(n\ell)| \ges |W((n+k)\ell)-W(0)  |-|W(n\ell))-W(0)| $$ 
for $k \in \N$, we see that
\begin{eqnarray}
q_n &\les& P(|W((n+k)\ell)-W(n\ell)| > (1+\ep)a^{n\ell} a_{n,k}   \mbox{ for } \forall k \in \N  ) \nonumber\\
& = &P(|W(k\ell) -W(0)| > (1+\ep)a_{n,k} \mbox{ for } \forall k \in \N  )\nonumber\\
&  \les& P(|W(k\ell) -W(0)| > b_{n,k} \mbox{ for } \forall k \in \N  ) =p_n=0.\nonumber
\end{eqnarray}
Thus we obtain that $ q_n = 0$ for every $n \in \Z_+$. Hence we have
$$P(\cup_{n=0}C_n)= \sum_{n=0}^{\infty}q_n =0.$$
Thus we established that
$$ \liminf_{n \to \infty}\frac{|W({n\ell}) -W(0)| }{a^{n\ell}g(n\ell) }\les 1+\ep.$$
Note that $\lim_{n \to \infty}a^ng(n)= \infty$ for $ g \in \G_1$. Since $\ep >0$ can be arbitrarily small, we have proved the lemma. \qed 

\begin{lem}  Let $g \in \G$. If $\sum_{n=0}^{\infty}P(|W(0)| \les g(n)) = \infty$, then $(3.6)$ holds.
\end{lem}
{\it Proof } Suppose that $\sum_{n=0}^{\infty}P(|W(0)| \les g(n)) = \infty$. We see from Lemma 3.5  that there exits $g^* \in \G_1$ such that $g(n) \les g^*(n)$ on $\Z_+$ and, for every $\delta > 0$, 
\begin{equation}
\sum_{n \in J}P(|W(0)| \les \delta g^*(n)) < \infty,
\end{equation}
where $J:=\{n \in\Z_+ : g(n) < g^*(n)\}.$
Since $\sum_{n=0}^{\infty}P(|W(0)| \les g^*(n)) = \infty$, we have by Lemma 3.6
\begin{equation}
\liminf_{n \to   \infty}\frac{|W(n)|}{a^ng^*(n)} \les 1 \mbox{ a.s.}.
\end{equation}
If the set $J$ is a finite set, then  (3.6) clearly holds by (3.11). If the set $J$ is an infinite set, then we can write $J$ as
\begin{equation}
J = \{ n_k : k \in \N, n_k \mbox{ is strictly increasing in } k \}.\nonumber
\end{equation}
We obtain from (3.10) as in the proof of Lemma 3.1 that
\begin{equation}
\liminf_{k \to   \infty}\frac{|W( n_k)|}{a^{n_k}g^*(n_k)} = \infty \mbox{ a.s.}.\nonumber
\end{equation}
Thus even in this case we have (3.6) by (3.11). \qed

\begin{thm} Let  $g \in \G.$ Then there is $C \in [0,\infty]$ such that
\begin{equation}
\liminf_{n \to   \infty}\frac{|W(n)|}{a^ng(n)} = C \mbox{ {\em  a.s.}}.
\end{equation}
The constant $C$ is determined by
$$ \sum_{n=0}^{\infty}P(|W(0)| \les \delta  g(n))  \left\{
\begin{array}{rl}  < \infty & \quad \mbox{ for $ 0 <\delta < C$}\\
 = \infty & \quad \mbox{ for $ \delta > C.$}
\end{array}\right. $$
\end{thm}
{\it Proof } The proof is clear from Lemmas 3.1 and 3.7.   \qed

\medskip

The following lemma is due to Remark 5.2 of Watanabe \cite{w02a}.
\begin{lem}  Let $b >a$. We have 
\begin{equation}
\lim_{n \to   \infty}\frac{|W(n)|}{b^n} = 0 \mbox{ {\em a.s.}}.\nonumber
\end{equation}
\end{lem}

\medskip
\begin{lem}  Let $a >b >1$. We have 
\begin{equation}
\lim_{n \to   \infty}\frac{|W(n)|}{b^n} = \infty \mbox{ {\em a.s.}}.\nonumber
\end{equation}
\end{lem}

{\it Proof }  We obtain from Lemma 3.4 that, for  every  $\delta >0$,
\begin{equation}
\sum_{n=0}^{\infty}P(|W(n)| \les \delta b^{n}) =\sum_{n=0}^{\infty}P(|W(0)| \les \delta b^{n}a^{-n})< \infty.\nonumber
\end{equation}
By virtue of lemma 3.1, we have proved the lemma.\qed 

\begin{thm} We have
\begin{equation}
\lim_{n \to   \infty}\frac{\log |W(n)|}{n} = \log a \mbox{ {\em  a.s.}}.\nonumber
\end{equation}
\end{thm}
{\it Proof } Proof  is due to Lemmas 3.8 and 3.9.  \qed

\medskip

The following lemma is due to Theorem 4.2 of Watanabe \cite{w02a}.

\begin{lem}  Let  $\{W(n)\}$ be an increasing   shift $a$-selfsimilar additive  sequence on $\R_+$. 
There exists $g \in \G$ satisfying $(1.4)$ if and only if  $P(W(0)-W(-1) =0)=0$. 
\end{lem}
\medskip

{\it Proof of Theorem 1.1. } 
First we prove (i).  We have
$$ \sum_{n=0}^{\infty}P(|W(0)| \les \delta) = \left\{
\begin{array}{rl} 0  & \quad \mbox{ for $ 0 <\delta < D$}\\
  \infty & \quad \mbox{ for $ \delta > D.$}
\end{array}\right. $$
Hence we established (i) from Theorem 3.1. Secondly we prove (ii). 
Suppose that $P(|W(0)| \les r)$ is  of dominated variation as $r \to 0+$ and (1.4) holds. We see from Theorem 3.1 that
\begin{equation}
 \sum_{n=0}^{\infty}P(|W(0)| \les 2g(n)) = \infty\nonumber
\end{equation}
and
\begin{equation}
 \sum_{n=0}^{\infty}P(|W(0)| \les 2^{-1}g(n)) < \infty.\nonumber
\end{equation}
This is a conradiction. Thus if (1.4) holds, then $P(|W(0)| \les r)$ is not  of dominated variation as $r \to 0+$. Next suppose that $P(|W(0)| \les r)$ is not of dominated variation as $r \to 0+$. 
Further, by virtue of (i), we can assume that $D=0$.  Then there exists a decreasing sequence $\{a_k\}_{k=1}^{\infty}$  such that $a_k >0$,  $a_{k+1} < a^{-1}a_k$ and 
\begin{equation}
 P(|W(0)| \les a^{-1}a_k) \les 2^{-k}P(|W(0)| \les a_k)
\end{equation}
for $k \ges 1$. We can choose an increasing sequence $\{b_k\}_{k=1}^{\infty}$ of integers such that $b_1=0$ and, for $k \ges 1$,
\begin{equation}
 2^{-1} \les P(|W(0)| \les a_k)(b_{k+1}-b_k) \les 1.
\end{equation}
 We define $g \in \G$ as
$ g(n):= a_k$ for $b_k \les n < b_{k+1}$ with $k \ges 1$. Then we have by (3.13) and (3.14) 
\begin{equation}
 \sum_{n=0}^{\infty}P(|W(0)| \les g(n)) \ges \sum_{k=1}^{\infty}2^{-1} = \infty\nonumber
\end{equation}
and 
\begin{equation}
 \sum_{n=0}^{\infty}P(|W(0)| \les a^{-1}g(n)) \les \sum_{k=1}^{\infty}2^{-k} < \infty.\nonumber
\end{equation}
Hence we see from Theorem 3.1  that (3.12)  holds for $C \in [a^{-1},1]$.  We have (1.4) by replacing $g(n)$  with $Cg(n)$. \qed
\medskip

{\it Proof of Proposition 1.1. } 
Let $W(n)= (W_j(n))_{j=1}^{d}$ and define an  increasing shift $a$-selfsimilar additive  sequence $\{Z(n)\}$ on $\R_+$ as
$$ Z(n):= \sum_{j=1}^d W_j(n).$$
Then we have
\begin{equation}
\frac{Z(n)}{\sqrt{d}}\les |W(n)| \les Z(n) 
\end{equation}
and 
\begin{equation}
\frac{Z(0)-Z(-1)}{\sqrt{d}}\les |W(0)-W(-1)| \les Z(0)-Z(-1). 
\end{equation}
Thus there exists $g \in \G$ satisfying (1.4) if and only if there exists $g \in \G$ such that
\begin{equation}
\liminf_{n \to   \infty}\frac{Z(n)}{a^ng(n)} = 1 \mbox{  a.s.}.\nonumber
\end{equation}
Therefore we see from Lemma 3.10 and  (3.16) that there exists $g \in \G$ satisfying (1.4) if and only if $P(Z(0)-Z(-1)=0)=0$, equivalently,
 $P(|W(0)-W(-1)| =0)=0$.  Thus the proposition has been proved. \qed

\medskip

\medskip

{\it Proof of Theorem 1.2. } Suppose that $P(|W(0)| \les r)$ is  of dominated variation as $r \to 0+$. Note that 
\begin{equation}
 \sum_{n=0}^{\infty}F(g(n)) = \infty\nonumber
\end{equation}
if and only if, for every $\delta >0$, 
\begin{equation}
 \sum_{n=0}^{\infty}F(\delta g(n)) = \infty.\nonumber
\end{equation}
Thus the theorem follows from Theorem 3.1. \qed

\medskip

{\it Proof of Proposition 1.2. } 
We continue to use $Z(n)$ as in the proof of Proposition 1.1. Suppose that $\lambda := P(|W(0)-W(-1)|=0 )= P(Z(0)-Z(-1)=0) >0.$
Define a Laplace transform $ L_Z(u)$ for $ u \ges 0$ as
$$L_Z(u):= E(\exp(-u(Z(0)-Z(-1))))$$ 
and a regularly varying function $K_{Z}(r)$ on $(0,1]$ with index $-\log \lambda /\log a$  as $r \to 0+$  as
$$ K_{Z}(r):=  r^{-\log \lambda/ \log a} \exp( \int_1^{r^{-1}} \frac{\log L_Z(u)-\log \lambda }{u \log a} du ). $$
Then we see from Proposition 4.1 of  Watanabe \cite{w02a} that $P(Z(0) \les r) \asymp K_{Z}(r)$ as $r \to 0+$.
We obtain from (3.16) that $K_{W}(r) \asymp K_{Z}(r)$ as $r \to 0+$.
Thus we see from (3.15) that, as $r \to 0+$,
\begin{equation}
P(|W(0)| \les r) \asymp P(Z(0) \les r)\asymp K_{Z}(r) \asymp K_{W}(r).\nonumber
\end{equation}
Therefore we have established Proposition 1.2 from  Theorem 1.2.\qed

\section{ Escape rates for strictly stable L\'evy  processes }
In this section, we show some laws of the iterated logarithm  for transient strictly stable L\'evy processes and complete the classical study of Taylor \cite{ta67} on the rates of escape for those processes.  Recurrent but not point recurrent cases are treated by Spitzer \cite{s} and Takeuchi and Watanabe \cite{tw}.

\begin{defn}
A L\'evy process $\{S(t), t \in \R_+\}$  on $\R^d$ is called a {\it strictly stable L\'evy process} if, for every $b>0$, there is $c>0$  such that
\begin{equation}
\{ S(bt)\cl t\ges0\}\eqd \{ cS(t)\cl t\ges0\}. \nonumber
\end{equation}
\end{defn}
 Let $\{S(t)\}$ be a strictly stable L\'evy  process on $\R^d$. We assume that the distribution of $S(1)$ is full. Then the distribution of $S(t)$ for $t >0$ is absolutely continuous with a density $p(t,x)$ of class $C^{\infty}$ in $x$. The index $\al$ is defined as $\al: = \log b/\log c$ for $b\ne1$ and satisfies  $0 < \alpha \les 2$.  

\begin{defn} The process $\{S(t)\}$ is called of {\it type I} if $p(1,0)>0$ and called of {\it type II} if $p(1,0)=0$.
 \end{defn}
\begin{defn} Let $\{S(t)\}$ be a strictly stable  L\'evy  process on a half space $\mathbb{H}^d:=\{ x : \langle z_0,x \rangle \ges 0 \mbox{ for some } z_0 \in \R^d\setminus\{0\} \}$ with index $0 < \alpha <1$. We define constants $C^*, C_* \in [0,\infty] $ by
$$ \sum_{n=1}^{\infty}P(|S(1)| \les \delta  (\log n)^{(\alpha-1)/\alpha}) \left\{
\begin{array}{rl}  < \infty & \quad \mbox{ for $ 0 <\delta < C^*$}\\
 = \infty & \quad \mbox{ for $ \delta > C^*,$}
\end{array}\right. $$
$$ \sum_{n=1}^{\infty}(\log n )^{1-\alpha}P(|S(1)| \les \delta  (\log n)^{(\alpha-1)/\alpha}) \left\{
\begin{array}{rl}  < \infty & \quad \mbox{ for $ 0 <\delta < C_*$}\\
 = \infty & \quad \mbox{ for $ \delta > C_*.$}
\end{array}\right. $$
\end{defn}
The process $\{S(t)\}$ is transient if and only if $d > \alpha$. Taylor \cite{ta67} proved that $\{S(t)\}$ is of type II if and only if $0 < \alpha < 1$ and the support of the distribution of $S(1)$ is included in a half space $\mathbb{H}^d$. Note that $\{S(t), t \in \R_+\}$ is a Sato process with exponent $1/\alpha$ and  $\{S(b^n), n \in \Z\}$ for $b > 1$ is a shift $b^{1/\alpha}$- selfsimilar additive sequence on $\R^d$. In the case  of type I, the escape rates of $\{S(b^t), t \in \R^1 \}$ and $\{S(b^n), n \in \Z\}$ are  different for $d > \alpha$ as the following proposition and its remark show.

\begin{prop}  Let $b > 1$ and $g(n) \in  \G $ with $n \in \Z_+$. If $\{S(t)\}$ is of type I and if
\begin{equation}
 \sum_{n=0}^{\infty}(g(n))^d = \infty \q ( \mbox{{\em resp.}} \q < \infty ),
\end{equation}
then
\begin{equation} 
\liminf_{n \to  \infty}\frac{|S(b^n)|}{b^{n/\alpha}g(n)} = 0 \q ( \mbox{{\em resp.}}  \q = \infty )\q \mbox{ {\em a.s..}}\nonumber
\end{equation}
In particular, if $\{S(t)\}$ is of type I, then for $\ep >0$ 
\begin{equation} 
\liminf_{n \to  \infty}\frac{|S(b^n)|}{b^{n/\alpha}n^{-1/d}(\log n )^{-(1+\ep)/d}} =  \infty  \mbox{ {\em a.s..}}\nonumber
\end{equation}

\end{prop}
{\it Proof } Note that, for every $x \in \R^d$, $P(|S(1)| \les r) \sim c r^d$ with some $c \in (0,\infty)$ as $r \to 0+$. Thus we can prove the proposition by using Theorem 1.2.  \qed 

\medskip                                                                                                                             
\begin{rem} Let $b > 1$.  Dvoretzky-Erd\"os \cite{de} for $\alpha =2$ and  Taylor \cite{ta67}  for $\alpha \ne 2$ showed the following.  Let   $\phi(t) $ be a positive decreasing function on $\R_+$. For $d > \alpha $, if $\{S(t)\}$ is of type I and if
\begin{equation}
 \int_{1}^{\infty}\phi(t)^{d-\alpha}dt = \infty \q ( \mbox{resp.} \q < \infty ),\nonumber
\end{equation}
then
\begin{equation} 
\liminf_{t \to \infty}\frac{|S(b^t)|}{b^{t/\alpha}\phi(t)} = 0 \q ( \mbox{resp.}  \q = \infty )\q \mbox{ a.s..}\nonumber
\end{equation}
In particular, if $\{S(t)\}$ is of type I, then for $d> \alpha$
\begin{equation} 
\liminf_{t \to \infty}\frac{|S(b^t)|}{b^{t/\alpha}t^{-1/(d-\alpha)}} = 0 \mbox{  a.s..}\nonumber
\end{equation}
\end{rem}

 The following lemma is a multi-dimensional extension of Lemmas 3.1 and 3.2  of Watanabe \cite{w96}

\begin{lem} Let $\{X(t), t \in \R_+\}$ be a Sato process on $\R_+^d$ with exponent $H $. Let $\phi(t)$ be a positive decreasing function on $\R_+$. Then there is $C \in [0,\infty]$ such that
\begin{equation}
\liminf_{t \to   \infty}\frac{|X(t)|}{t^H\phi(t)} = C \mbox{ {\em a.s.}}.
\end{equation}
The constant $C$ is determined by
$$ \int_{1}^{\infty}P(|X(1)| \les \delta  \phi(t))t^{-1}dt  \left\{
\begin{array}{rl}  < \infty & \quad \mbox{ for $ 0 <\delta < C$}\\
 = \infty & \quad \mbox{ for $ \delta > C.$}
\end{array}\right. $$
\end{lem}
{\it Proof } Let $b >1$.    Note that $\phi(b^n) \in \G$ and $\{X(b^n)\}$ is a shift $b^H$-selfsimilar additive sequence on $\R_+^d$. We define $C_0 \in [0,\infty]$  by 
$$ \int_{1}^{\infty}P(|X(1)| \les \delta  \phi(t))t^{-1}dt  \left\{
\begin{array}{rl}  < \infty & \quad \mbox{ for $ 0 <\delta < C_0$}\\
 = \infty & \quad \mbox{ for $ \delta > C_0.$}
\end{array}\right. $$
Since $\phi(t)$ is decreasing, we see that 
$$ \sum_{n=0}^{\infty}P(|X(1)| \les \delta  \phi(b^n))  \left\{
\begin{array}{rl}  < \infty & \quad \mbox{ for $ 0 <\delta < C_0$}\\
 = \infty & \quad \mbox{ for $ \delta > C_0.$}
\end{array}\right. $$
Thus we obtain from Theorem 3.1 that
\begin{equation}
\liminf_{n \to   \infty}\frac{|X(b^n)|}{b^{nH}\phi(b^n)} = C_0 \mbox{ a.s.}.\nonumber
\end{equation}
Since $|X(t)|$ is increasing and $\phi(t)$ is decreasing, we have  for $b^n < t \les b^{n+1}$ 
\begin{equation}
\frac{|X(b^n)|}{b^{(n+1)H} \phi(b^n)}   \les \frac{|X(t)|}{t^H\phi(t)} \les \frac{|X(b^{n+1})|}{b^{nH}\phi(b^{n+1})}   \mbox{ a.s.}\nonumber
\end{equation}
and thereby
\begin{equation}
b^{-H} C_0  \les \liminf_{t \to   \infty}\frac{|X(t)|}{t^H\phi(t)} \les b^{H} C_0  \mbox{ a.s.}.\nonumber
\end{equation}
Letting $b \to 1$,  we have (4.2) with $C=C_0$. \qed 

\begin{lem} Let  $\{S(t)\}$ be a strictly stable  L\'evy  process  on $\R_+^d$ with index $0 < \alpha <1$.   Then we have as $r \to 0+$
\begin{equation}
 -\log P(|S(1)| \les r ) \asymp r^{\alpha/(\alpha-1)}.
\end{equation}
\end{lem}

{\it Proof } Let $S(t)=(S_j(t))_{j=1}^d$. Define a strictly stable  L\'evy  process $\{S^*(t)\}$ on $\R_+$  with index $0 < \alpha <1$ as $S^*(t):= \sum_{j=1}^dS_j(t)$. We see from Remark 14.18  of Sato \cite{s99} that there exists $c >0$ such that as $r \to 0+$
\begin{equation}
 -\log P(|S^*(1)| \les r ) \sim c r^{\alpha/(\alpha-1)}.\nonumber
\end{equation}
Thus since 
\begin{equation}
\frac{|S^*(1)|}{\sqrt{d}}\les |S(1)| \les |S^*(1)|,\nonumber
\end{equation}
we have (4.3) as $r \to 0+$. \qed

\begin{lem} Let  $\{S(t)\}$ be a strictly stable  L\'evy  process  on a half space $\mathbb{H}^d$ with index $0 < \alpha <1$.   Then   we have $(4.3)$ as $r \to 0+$.
\end{lem}
{\it Proof } Without loss of generality, we can assume that $z_0={}^t(1,0 \dots 0) \in \R^d$ and $\mathbb{H}^d= \R_+\times \R^{d-1}.$ We define a set $\D$ of disjoint orthants as 
\begin{equation}
\D:= \{D : D=\R_+\times E_1 \times \cdots \times E_{d-1},  E_k= [0,\infty)\mbox{ or } (-\infty,0) \mbox{ for } 1 \les k \les d-1 \}.\nonumber
\end{equation}
Let $\{D_j\}_{j=1}^{N}$ be the totality of $D_j \in \D$ such that $\nu(D_j) \ne 0$. Let $\nu$ be the L\'evy measure of the process $\{S(t)\}$ and define $\nu_j$ for $1 \les j \les N$ as $\nu_j(dx):= 1_{D_j}(x)\nu(dx)$. Let $\{X_j(t)\}$ for $1 \les j \les N$ be independent strictly stable L\'evy processes with L\'evy measures $\nu_j$ and index $0 < \alpha <1$. Then we have 
\begin{equation}
\{ S(t)\cl t\ges0\}\eqd \{ \sum_{j=1}^{N}X_j(t)\cl t\ges0\}. \nonumber
\end{equation}
Define a strictly stable  L\'evy  process $\{S_*(t)\}$ on $\R_+$  with index $0 < \alpha <1$ as $S_*(t):= \langle z_0,S(t) \rangle$. 
Since $|S_*(1)| \les |S(1)| \les \sum_{j=1}^{N}|X_j(1)|$, we have
\begin{equation}
\prod_{j=1}^NP(|X_j(1)|\les N^{-1}r) \les P(|S(1)| \les r )\les P(|S_*(1)| \les r ).\nonumber
\end{equation}
Thus from Lemmas 4.2 we obtain (4.3)  as $r \to 0+$. \qed

\begin{lem} Let  $\{S(t), t \in \R_+\}$ be a strictly stable  L\'evy  process  on a half space $\mathbb{H}^d$ with index $0 < \alpha <1$.  Then we have
\begin{equation}
0 < C_* \les  C^* < \infty.\nonumber
\end{equation}
\end{lem}
{\it Proof } The proof is clear from Lemma 4.3. \qed

\medskip

\begin{thm} Let $\{S(t)\}$ be a strictly stable  L\'evy  process  on $\R_+^d$ with index $0 < \alpha <1$. Then  we have
\begin{equation}
\liminf_{t \to   \infty}\frac{|S(t)|}{t^{1/\alpha}(\log \log t)^{(\alpha-1)/\alpha}} = C^* \mbox{ {\em a.s.}}.
\end{equation}
\end{thm}

{\it Proof } Let $\phi(t)= (\log \log t)^{(\alpha-1)/\alpha}$. Then we obtain (4.4) from (4.2) of Lemma 4.1. The  constant $C$ in (4.2) is determined by
$$ \int_{e}^{\infty}P(|S(1)| \les \delta  (\log \log t)^{(\alpha-1)/\alpha})t^{-1}dt  \left\{
\begin{array}{rl}  < \infty & \quad \mbox{ for $ 0 <\delta < C$}\\
 = \infty & \quad \mbox{ for $ \delta > C.$}
\end{array}\right. $$
That is, $C = C^*$.  \qed

\medskip

The following lemma is a modification of Theorem 1.1 of Khoshnevisan \cite{k97}. 
\begin{lem}
Let $\{Z(t)\}$ be a L\'evy process on $\R^d$.  For any $0 < b < c$, $\gamma >0$, and $\ep >0$, we have
\begin{equation}
P(|Z(t)| \les \gamma \mbox{ for some } b \les t \les c) \les \frac{\int_b^{2c-b}P(|Z(t)| \les (1+\ep)\gamma)dt}{\int_0^{c-b}P(|Z(t)| \les \ep \gamma)dt}\nonumber
\end{equation}   

\end{lem}
{\it Proof } Define $T:=\inf \{t \ges b : |Z(t)| \les \gamma \}$. Then we see that 
\begin{equation} 
\begin{split}
&\int_b^{2c-b}P(|Z(t)| \les (1+\ep)\gamma)dt\nonumber \\
&\ges P(T \les c)\int_b^{2c-b}P(|Z(t)| \les (1+\ep)\gamma| T \les c)dt\\
&\ges P(T \les c)\inf_{|x| \les \gamma}\int_0^{c-b}P(|x+Z(t)| \les (1+\ep)\gamma)dt\\
&\ges P(T \les c)\int_0^{c-b}P(|Z(t)| \les \ep\gamma)dt\\
\end{split}
\end{equation}
Thus we have proved the lemma. \qed

\medskip

The following is a multi-dimensional extension of a result of  Fristedt \cite{f64}.

\begin{thm} Let $\{S(t)\}$ be a strictly stable  L\'evy  process  on  a half space $\mathbb{H}^d$ with index $0 < \alpha <1$. Then $0 < C_* \les  C^* < \infty$ and there is $C \in [C_*, C^* ]$ such that
\begin{equation}
\liminf_{t \to   \infty}\frac{|X(t)|}{t^{1/\alpha}(\log \log t)^{(\alpha-1)/\alpha}} = C \mbox{ {\em a.s.}}.
\end{equation}
\end{thm}

{\it Proof } Lemma 4.4 says that $0 < C_* \les  C^* < \infty$. By virtue of Kolmogorov's 0-1 law, we have (4.5) for $C \in [0,\infty]$. Let $W(n) := S(e^n)$ for $n \in \Z$ and 
$$g(n):= (\log \log e^n)^{(\alpha-1)/\alpha}=(\log n)^{(\alpha-1)/\alpha}$$
for $n \ges 2$.  Then we see from the definition of $C^*$ and Theorem 3.1  that
\begin{equation}
\liminf_{n \to   \infty}\frac{|S(e^n)|}{e^{n/\alpha}(\log \log e^n)^{(\alpha-1)/\alpha}} = C^* \mbox{  a.s.}\nonumber
\end{equation}
 Thus we obtain that $C \les C^*.$ Fix a small $\ep \in (0,1)$. Then we choose $b >1$ satisfying  
\begin{equation}
1-\ep < b^{1/\alpha}(1-\ep^2) < 1.\nonumber
\end{equation}
Let 
$$\gamma_n:= C_*(1-\ep) b^{(n+1)/\alpha}(\log (n+1))^{(\alpha-1)/\alpha}$$ 
for $n \in \N$ and put
$$c_1:=\int_0^{\infty}P(|S(1)| \les u^{-1/\alpha}\ep )du.$$
  We see from Lemma 4.3 that $c_1 < \infty$. Thus there exists $c_2 >0$ such that, for $n \in \N$,
\begin{equation} 
\begin{split}
&I_1(n):=\int_0^{b^{n+1}-b^n}P(|S(t)| \les \ep \gamma_n)dt  \\
=&\int_0^{b^{n+1}-b^n}P(|S(1)| \les t^{-1/\alpha}\ep \gamma_n)dt  \\
=&\gamma_n^{\alpha}\int_0^{(b^{n+1}-b^n)\gamma_n^{-\alpha}}P(|S(1)| \les u^{-1/\alpha}\ep )du.  \\
\ges & c_2 b^{n+1}(\log (n+1))^{\alpha-1}. \\
\end{split}
\end{equation}
On the other hand, we have for $n \in \N$
\begin{equation} 
\begin{split}
&I_2(n):=\int_{b^n}^{2b^{n+1}-b^n}P(|S(t)| \les (1+\ep) \gamma_n)dt  \\
=&\int_{b^n}^{2b^{n+1}-b^n}P(|S(1)| \les t^{-1/\alpha}(1+\ep) \gamma_n)dt  \\
\les&2(b^{n+1}-b^n)P(|S(1)| \les C_* b^{1/\alpha}(1-\ep^2)(\log (n+1))^{(\alpha-1)/\alpha} ).  \\
\end{split}
\end{equation}
Thus we obtain from Lemma 4.5, (4.6), and (4.7) that there exists $c_3 >0$ such that
\begin{equation}
\begin{split}
&\sum_{n=1}^{\infty}P(|S(t)| \les \gamma_n \mbox{ for some } b^n \les t \les b^{n+1})\nonumber  \\
\les &\sum_{n=1}^{\infty}\frac{I_2(n)}{I_1(n)}\\
\les &\sum_{n=1}^{\infty}c_3(\log (n+1))^{1-\alpha}P(|S(1)| \les  C_* b^{1/\alpha}(1-\ep^2)(\log (n+1))^{(\alpha-1)/\alpha} ) < \infty.
\end{split}
\end{equation} 
Thus by the Borel-Cantelli lemma, we have 
\begin{equation}
\liminf_{n \to   \infty}\inf_{b^n \les t \les  b^{n+1}}\frac{|S(t)|}{b^{(n+1)/\alpha}(\log (n+1))^{(\alpha-1)/\alpha} } \ges C_*(1-\ep) \mbox{  a.s.}.\nonumber
\end{equation}
Since $\ep >0$ is arbitrarily small, we have $C \ges C_*$ and hence $ C_*\les C \les C^* .$ \qed 

\section{LIL for independent  Brownian motions }
In this section, we show some laws of the iterated logarithm (LIL in short) for independent Brownian motions with an explicit constant.  
Define the norm $|x|^*:= \sup_{1 \les j \les d}|x_j|$ for $x =(x_j)_{j=1}^d \in \R^d$. We see from the proofs of previous results that the use of the norm $|x|^*$ instead of the norm $|x|$ does not harm the validity of  the following lemmas. The first one is due to Lemma 4.1.

\begin{lem} Let $\{X(t), t \in \R_+\}$ be a Sato process on $\R_+^d$ with exponent $H $. Let $\phi(t)$ be a positive decreasing function on $\R_+$. Then there is $C \in [0,\infty]$ such that
\begin{equation}
\liminf_{t \to   \infty}\frac{|X(t)|^*}{t^H\phi(t)} = C \mbox{ {\em a.s.}}.
\end{equation}
The constant $C$ is determined by
$$ \int_{1}^{\infty}P(|X(1)|^* \les \delta  \phi(t)))t^{-1}dt  \left\{
\begin{array}{rl}  < \infty & \quad \mbox{ for $ 0 <\delta < C$}\\
 = \infty & \quad \mbox{ for $ \delta > C.$}
\end{array}\right. $$
\end{lem}
\medskip 

\medskip

\begin{lem} Let $\{X(t), t \in \R_+\}$ be a Sato process on $\R_+^d$ with exponent $H $. 
We have
\begin{equation} 
\lim_{t \to  \infty}\frac{\log (|X(t)|^*)}{\log t} =  H \mbox{ {\em a.s..}}
\end{equation}

\end{lem}

{\it Proof } Let $b >1$. Since $|X(t)|^*$ is increasing, we have for $b^n < t \les b^{n+1}$ 
\begin{equation}
\frac{\log (|X(b^{n})|^*)}{\log b^{n+1}}  \les \frac{\log (|X(t)|^*)}{\log t}  \les \frac{\log (|X(b^{n+1})|^*)}{\log b^n}   \mbox{ a.s.}.\nonumber
\end{equation}
We see from Theorem 3.2 that
\begin{equation}
\lim_{n \to \infty}\frac{\log (|X(b^{n})|^*)}{n}  = H \log b   \mbox{ a.s.}.\nonumber
\end{equation}
Thus we obtain (5.2). \qed

\medskip

 Let $\{B_j(t)\}$ be independent  Brownian motions on $\R^d$ starting at the origin for $1 \les j \les N$ with $N \in \N$. For $r \ges 0$, define the first hitting time $\{T_j(r)\}$ to the sphere $S^{d-1}(r): =\{x : |x| =r\}$ for $d \ges 1$ and the last exit time $\{L_j(r)\}$  from  $S^{d-1}(r)$ for $d \ges 3$ as
\begin{equation} 
T_j(r):= \inf \{t \ges 0 : |B_j(t)| = r \}\nonumber
\end{equation} 
and 
\begin{equation} 
L_j(r):= \sup \{t \ges 0 : |B_j(t)| = r \}.\nonumber
\end{equation} 
Let $T(t):= (T_j(t))_{j=1}^N$ and $L(t):= (L_j(t))_{j=1}^N$. Then, by virtue of  Getoor \cite{g79},  $\{T(t)\}$ and $\{L(t)\}$  are Sato processes on $\R_+^N$ with $H=2$.

\begin{rem} We have by Lemma 5.2 
\begin{equation} 
\lim_{r \to  \infty}\frac{\log (\sup_{1 \les j \les N}T_j(r))}{\log r} =  2 \mbox{  a.s.}
\end{equation}
and, for $d \ges 3$,
\begin{equation} 
\lim_{r \to  \infty}\frac{\log (\sup_{1 \les j \les N}L_j(r))}{\log r} =  2 \mbox{  a.s..}\nonumber
\end{equation}
\end{rem}

\begin{lem}  Let $c_1 , c_2$ be positive constants. We have as $r \to 0+$
\begin{equation} 
 P(T_1(1) \les r) \sim c_1 \exp(-1/(2r))r^{-d/2 +1} 
\end{equation}
and
\begin{equation} 
 P( L_1(1)  \les r) \sim c_2 \exp(-1/(2r))r^{-d/2 +2}.
\end{equation}

\end{lem}
{\it Proof } The relation (5.4) is due to Theorem 2.1 of  Gruet and Shi \cite{gs} and (5.5) is due to Theorem of Getoor \cite{g79}.
\qed

\begin{thm}  We have
\begin{equation} 
\liminf_{r \to  \infty}\frac{\sup_{1 \les j \les N}T_j(r)}{r^2(\log\log r)^{-1}} = \frac{N}{2} \mbox{ {\em a.s.}}
\end{equation}
 and, for $d \ges 3$,
\begin{equation} 
\liminf_{r \to  \infty}\frac{\sup_{1 \les j \les N}L_j(r)}{r^2(\log\log r)^{-1}} = \frac{N}{2} \mbox{ {\em a.s..}}
\end{equation}
\end{thm}
{\it Proof.}  We use the norm $|x|^*:= \sup_{1 \les j \les N}|x_j|$ for $x =(x_j)_{j=1}^N \in \R_+^N$.   We have by (5.4) of Lemma 5.3
$$ P(|T(1)|^* \les r) = (P(T_1(1) \les r))^N \sim c_1^N \exp(-N/(2r))r^{(-d/2 +1)N}.$$
Thus we see that
$$ \int_{e}^{\infty}P(|T(1)|^* \les \delta  (\log \log t)^{-1})t^{-1}dt  \left\{
\begin{array}{rl}  < \infty & \quad \mbox{ for $ 0 <\delta < N/2$}\\
 = \infty & \quad \mbox{ for $ \delta > N/2.$}
\end{array}\right. $$
Hence we obtain (5.6) from (5.1) of Lemma 5.1.   
In the same way, we have (5.7) by (5.5) of Lemma 5.3.  \qed

\begin{thm}  We have
\begin{equation} 
\limsup_{t \to  \infty}\inf_{1 \les j \les N}\frac{\sup_{0 \les s \les t}|B_j(s)|}{\sqrt{t \log\log t}} = \sqrt{\frac{2}{N}} \mbox{ {\em a.s.}}
\end{equation}
 and, for $d \ges 3$,
\begin{equation} 
\limsup_{t \to  \infty}\inf_{1 \les j \les N}\frac{\inf_{t \les s}|B_j(s)|}{\sqrt{t \log\log t}} = \sqrt{\frac{2}{N}} \mbox{ {\em a.s..}}
\end{equation}

\end{thm}

{\it Proof } We prove only (5.8). The proof of (5.9) is analogous by using (5.7) of Theorem 5.1 and is omitted. Thanks to Kolmogorov's 0-1 law, we find that there is a constant $C \in [0,\infty]$ such that
\begin{equation} 
\limsup_{t \to  \infty}\inf_{1 \les j \les N}\frac{\sup_{0 \les s \les t}|B_j(s)|}{\sqrt{t \log\log t}} =C \mbox{  a.s..}
\end{equation}
We see from (5.6) of Theorem 5.1  that there exists a random sequence $\{r_n\}=\{r_n(\omega)\}$ such that $r_n \to \infty$ as $n \to \infty$ and that
\begin{equation} 
\lim_{n \to  \infty}\frac{\sup_{1 \les j \les N}T_j(r_n)}{r_n^2(\log\log r_n)^{-1}} = \frac{N}{2} \mbox{  a.s..}
\end{equation}
Denote $t_n :=\sup_{1 \les j \les N}T_j(r_n)$. Then we have 
\begin{equation} 
r_n =\inf_{1 \les j \les N}\sup_{0 \les s \les t_n}|B_j(s)|.\nonumber
\end{equation}
We find from (5.3) of  Remark 5.1 that 
\begin{equation} 
\lim_{n \to  \infty}\frac{\log t_n}{\log r_n} =  2 \mbox{  a.s..}\nonumber
\end{equation}
Thus we have by (5.11)
\begin{equation} 
\lim_{n \to  \infty}\frac{t_n}{r_n^2(\log\log t_n)^{-1}} = \frac{N}{2} \mbox{  a.s.}\nonumber
\end{equation}
and hence $C \ges \sqrt{\frac{2}{N}}.$
Suppose that there exists a random sequence $\{s_n\}=\{s_n(\omega)\}$ such that $s_n \to \infty$ as $n \to \infty$ and that
\begin{equation} 
C=\lim_{n \to  \infty}\inf_{1 \les j \les N}\frac{\sup_{0 \les s \les s_n}|B_j(s)|}{(s_n \log\log s_n)^{1/2}} > \sqrt{\frac{2}{N}} \mbox{  a.s..}
\end{equation}
Denote $u_n:=\inf_{1 \les j \les N}\sup_{0 \les s \les s_n}|B_j(s)|.$ Then we have 
\begin{equation} 
s_n \ges \sup_{1 \les j \les N}T_j(u_n).
\end{equation}
We see from (5.3) of Remark 5.1  that 
\begin{equation} 
\lim_{n \to  \infty}\frac{\log (\sup_{1 \les j \les N}T_j(u_n))}{\log u_n} =  2 \mbox{  a.s..}\nonumber
\end{equation}
Thus we have by (5.12) and (5.13)
\begin{equation} 
\liminf_{n \to  \infty}\frac{\sup_{1 \les j \les N}T_j(u_n)}{u_n^2(\log\log u_n)^{-1}} \les \frac{1}{C^2} < \frac{N}{2} \mbox{  a.s..}\nonumber
\end{equation}
This is a contradiction. Hence we have established that $C = \sqrt{\frac{2}{N}}$ in (5.10). \qed
\begin{rem}  

The symbol  ``$\sup_{0\le s \le t}|B(s)|$'' in (5.8) of Theorem 5.2 can be replaced by $|B(t)|$. For $N=1$,  (5.8) is due to Khintchine \cite{kh33} and  (5.9) is  to Khoshnevisan et al.\ \cite{kl94}.

\end{rem}

\end{document}